\documentclass[a4paper,oneside]{article}
\usepackage{amssymb}
\usepackage{amsfonts}
\usepackage{amsmath}

\input{tcilatex}

\begin{document}

\title{On the Borel transgression in the fibration $G\rightarrow G/T$}
\author{Haibao Duan\thanks{%
The author's research is supported by NSFC 11131008; 11661131004} \\
Institute of Mathematics, Chinese Academy of Sciences;\\
School of Mathematical Sciences, \\
University of the Chinese Academy of Sciences\\
dhb@math.ac.cn}
\maketitle

\begin{abstract}
Let $G$ be a semisimple Lie group with a maximal torus $T$. We present an
explicit formula for the Borel transgression $\tau :H^{1}(T)\rightarrow
H^{2}(G/T)$ of the fibration $G\rightarrow G/T$. This formula corrects an
error in the paper \cite{K}, and has been applied to construct the integral
cohomology rings of compact Lie groups in the sequel works \cite{D,DZ2}.

\begin{description}
\item[2010 Mathematical Subject Classification: ] 55T10, 57T10

\item[Key words and phrases:] Lie groups; Cohomology, Leray--Serre spectral
sequence
\end{description}
\end{abstract}

\section{Introduction}

A Lie group is called \textsl{semisimple} if its center is finite; is called 
\textsl{adjoint }if its center is trivial. In this paper the Lie groups $G$
under consideration are compact, connected and semisimple. The homology and
cohomology are over the ring of integers, unless otherwise stated.

For a Lie group $G$ with a maximal torus $T$ let $\pi :G\rightarrow G/T$ be
the quotient fibration. Consider the diagram with top row the cohomology
exact sequence of the pair $(G,T)$

\begin{quote}
\begin{tabular}{lllll}
$0\rightarrow H^{1}(G)\overset{i^{\ast }}{\rightarrow }$ & $H^{1}(T)$ & $%
\overset{\delta }{\rightarrow }$ & $H^{2}(G,T)$ & $\overset{j^{\ast }}{%
\rightarrow }H^{2}(G)\rightarrow \cdots $ \\ 
&  & $\searrow \tau $ & $\quad \cong \uparrow \pi ^{\ast }$ &  \\ 
&  &  & $H^{2}(G/T)$ & 
\end{tabular}
\end{quote}

\noindent where, since the pair $(G,T)$ is $1$--connected, the induced map $%
\pi ^{\ast }$ is an isomorphism. The \textsl{Borel} \textsl{transgression }%
\cite[p.185]{Mc}\textbf{\ }in the fibration $\pi $ is the composition

\begin{quote}
$\tau =$ $(\pi ^{\ast })^{-1}\circ \delta :H^{1}(T)\rightarrow H^{2}(G/T)$.
\end{quote}

\noindent It is irrelevant to the choice of a maximal torus on $G$ since, if 
$T^{\prime }$ is another maximal torus then the relation $T^{\prime
}=gTg^{-1}$ holds on $G$ for some $g\in $ $G$.

The transgression $\tau $ is an essential ingredient of the Leray--Serre
spectral sequence $\left\{ E_{r}^{\ast ,\ast }(G;\mathcal{R}),d_{r}\mid
r\geq 2\right\} $ of the fibration $\pi $, where the coefficient ring $%
\mathcal{R}$ is either the ring $\mathbb{Z}$ of integers, the field $\mathbb{%
R}$ of reals, or the cyclic ring $\mathbb{Z}_{p}$ of finite order $p$.
Precisely, as the quotient manifold $G/T$ is always $1$--connected, the
Leray--Serre Theorem states that

\begin{enumerate}
\item[(1.1)] $E_{2}^{\ast ,\ast }(G;\mathcal{R})=H^{\ast }(G/T)\otimes
H^{\ast }(T;\mathcal{R})$;

\item[(1.2)] the differential $d_{2}:E_{2}^{\ast ,\ast }\rightarrow
E_{2}^{\ast ,\ast }$ is determined by $\tau $ as
\end{enumerate}

\begin{quote}
$\qquad d_{2}(x\otimes t)=(x\cup \tau (t))\otimes 1$,
\end{quote}

\noindent where $x\in H^{\ast }(G/T)$, $t\in H^{1}(T)$. Furthermore,
regarding both $E_{3}^{\ast ,\ast }(G;\mathcal{R})$ and $H^{\ast }(G;%
\mathcal{R})$ as graded groups, Leray and Reeder \cite{L,Re} have shown that

\begin{quote}
$E_{3}^{\ast ,\ast }(G;\mathbb{R})=H^{\ast }(G;\mathbb{R})$;
\end{quote}

\noindent Ka\v{c} \cite{K} claimed that, if $p$ is a prime, then

\begin{quote}
$E_{3}^{\ast ,\ast }(G;\mathbb{Z}_{p})=H^{\ast }(G;\mathbb{Z}_{p})$.
\end{quote}

\noindent Marlin \cite{M} conjectured that, if the group $G$ is $1$%
--connected, then

\begin{quote}
$E_{3}^{\ast ,\ast }(G;\mathbb{Z})=H^{\ast }(G)$.
\end{quote}

\noindent Conceivably, an explicit formula for the transgression $\tau $ is
requested by the spectral sequence approach to the cohomology theory of Lie
groups.

Our main result is Theorem 2.5, where with respect to explicitly constructed
bases on $H^{2}(G/T)$ and $H^{1}(T)$ a formula for $\tau $ is obtained. This
formula corrects an error concerning the differential $d_{2}$ on $%
E_{2}^{\ast ,\ast }$ occurring in \cite{K}. It has also been applied in our
sequel works \cite{D,DZ2} to construct the integral cohomology rings of
compact Lie groups, see Remark 3.4.

\section{A formula for the transgression $\protect\tau $}

For a semisimple Lie group $G$ with a maximal torus $T$ the tangent space $%
L(G)$ (resp. $L(T)$) to $G$ at the group unit $e\in G$ (resp. to $T$ at $%
e\in T$) is also known as the \textsl{Lie algebra} (resp. the \textsl{Cartan
subalgebra}) of $G$. The exponential map $\exp $ of $G$ at the unit $e$
builds up the commutative diagram

\begin{quote}
\begin{tabular}{lll}
$\quad L(T)$ & $\rightarrow $ & $L(G)$ \\ 
$\exp \downarrow $ &  & $\downarrow \exp $ \\ 
$\quad T$ & $\rightarrow $ & $G$%
\end{tabular}%
,
\end{quote}

\noindent where the horizontal maps are the obvious inclusions. Equip the
real vector space $L(G)$ with an inner product $(,)$ so that the adjoint
representation acts as isometries on $L(G)$, and assume that $n=\dim T$.

The exponential map $\exp :L(T)\rightarrow T$ defines a set $\mathcal{S}%
(G)=\{L_{1},\cdots ,L_{m}\}$ of $m=\frac{1}{2}(\dim G-n)$ hyperplanes on $%
L(T)$, namely, the set of\textsl{\ singular hyperplanes }through the origin
in $L(T)$ \cite[p.168]{BT}. The map $\exp $ carries the normal line $l_{k}$
to $L_{k}$ through the origin $0\in L(T)$ onto a circle subgroup on $T$. Let 
$\pm \alpha _{k}\in l_{k}$ be the non--zero vectors with minimal length so
that

\begin{quote}
$\exp (\pm \alpha _{k})=e$, $1\leq k\leq m$.
\end{quote}

\noindent The subset $\Phi =\{\pm \alpha _{k}\in L(T)\mid 1\leq k\leq m\}$
will be called \textsl{the root system of }$G$.

In addition, the planes in $\mathcal{S}(G)$ divide $L(T)$ into finitely many
convex regions, called the \textsl{Weyl chambers} of $G$. Fix a regular
point $x_{0}\in L(T)$, and let $\mathcal{F}(x_{0})$ be the closure of the
Weyl chamber containing $x_{0}$. Assume that $L(x_{0})=\{L_{1},\cdots
,L_{n}\}\subset \mathcal{S}(G)$ is the subset consisting of the walls of $%
\mathcal{F}(x_{0})$, and let $\alpha _{i}\in \Phi $ be the root normal to $%
L_{i}$ and pointing toward $x_{0}$.

\bigskip

\noindent \textbf{Definition 2.1.} The subset $\Delta =\{\alpha _{1},\cdots
,\alpha _{n}\}\subset L(T)$ will be called \textsl{the set of simple roots}
of $G$ relative to the Weyl chamber $\mathcal{F}(x_{0})$ \cite[p.49]{H}.

The \textsl{Cartan matrix} of $G$ is the $n$ by $n$ matrix $%
A=(b_{ij})_{n\times n}$ defined by $b_{ij}=2(a_{i},\alpha _{j})/(\alpha
_{j},\alpha _{j})$.$\square $

\bigskip

\noindent \textbf{Remark 2.2. }Our geometric approach to the simple roots,
as well as the Cartan matrix of $G$,\textbf{\ }is \textsl{dual} to those
that are commonly used in literatures, e.g. \cite{B,H}. This allows us to
perform subsequent construction and calculation coherently on the Euclidean
space $L(T)$ without referring to its dual space $L(T)^{\ast }$.$\square $

\bigskip

On the Euclidean space $L(T)$ there are three distinguished lattices.
Firstly, the set $\Delta =\{\alpha _{1},\cdots ,\alpha _{n}\}$ of simple
roots generates over the integers \textsl{the} \textsl{root lattice} $%
\Lambda _{r}$ of $G$. Next, the pre--image of the exponential map $\exp
:L(T)\rightarrow T$ at the group unit $e\in T$ gives rise to the \textsl{%
unit lattice }$\Lambda _{e}:=\exp ^{-1}(e)$\textsl{\ }of $G$. Thirdly, using
simple roots one defines the set $\Omega =\{\phi _{1},\cdots ,\phi
_{n}\}\subset L(T)$ of \textsl{fundamental dominant weights} of $G$ by the
formula

\begin{quote}
$2(\phi _{i},\alpha _{j})/(\alpha _{j},\alpha _{j})=\delta _{i,j}$
\end{quote}

\noindent that generates over the integers \textsl{the weight lattice} $%
\Lambda _{\omega }$ of $G$. It is known that (see \cite[(3.4)]{DL}):

\bigskip

\noindent \textbf{Lemma 2.3. }\textsl{On the Euclidean space }$L(T)$ \textsl{%
one has }

\begin{quote}
$\Lambda _{r}\subseteq \Lambda _{e}\subseteq \Lambda _{\omega }$\textsl{. }
\end{quote}

\noindent \textsl{In addition}

\textsl{i) the group }$G$\textsl{\ is }$1$\textsl{--connected if and only if 
}$\Lambda _{r}=\Lambda _{e}$\textsl{;}

\textsl{ii) the group }$G$\textsl{\ is adjoint if and only if }$\Lambda
_{e}=\Lambda _{\omega }$\textsl{;}

\textsl{iii) the basis }$\Delta $\textsl{\ on }$\Lambda _{r}$\textsl{\ can
be expressed by the basis }$\Omega $\textsl{\ on }$\Lambda _{\omega }$ 
\textsl{by}

\begin{quote}
$\qquad \left( 
\begin{tabular}{l}
$\alpha _{1}$ \\ 
$\vdots $ \\ 
$\alpha _{n}$%
\end{tabular}%
\right) =A\left( 
\begin{tabular}{l}
$\phi _{1}$ \\ 
$\vdots $ \\ 
$\phi _{n}$%
\end{tabular}%
\right) $,
\end{quote}

\noindent \textsl{where }$A$\textsl{\ is the Cartan matrix of }$G$\textsl{.}$%
\square $

\bigskip

Granted with the notion introduced above we turn to the construction of a
canonical basis of the second cohomology $H^{2}(G/T)$ of the base manifold $%
G/T$. For a simple root $\alpha \in \Delta $ let $K(\alpha )\subset G$ be
the subgroup with Lie algebra $l_{\alpha }\oplus L_{\alpha }$ (\cite[p.238,
Exercise 6]{BT}), where $l_{\alpha }\subset L(T)$ is the $1$--dimensional
subspace spanned by $\alpha $, and $L_{\alpha }\subset L(G)$ is the root
space (viewed as an oriented real $2$--plane) belonging to the root $\alpha $
(\cite[p.35]{H}). Then the circle subgroup $S^{1}=\exp (l_{\alpha })$ is a
maximal torus on $K(\alpha )$, while the quotient manifold $K_{\alpha
}/S^{1} $ is diffeomorphic to the $2$--dimensional sphere $S^{2}$. Moreover,
the inclusion $(K_{\alpha },S^{1})\subset (G,T)$ of subgroups induces an
embedding

\begin{enumerate}
\item[(2.1)] $s_{\alpha }:S^{2}=K_{\alpha }/S^{1}\rightarrow G/T$
\end{enumerate}

\noindent whose image is known as the \textsl{Schubert variety} associated
to the root $\alpha $ \cite{DZ1}. By the basis theorem of Chevalley \cite{Ch}
the maps $s_{\alpha }$ with $\alpha \in \Delta $ represent a basis of the
second homology $H_{2}(G/T)$. As a result, if one lets $\omega _{i}\in
H^{2}(G/T)$\ be the Kronecker dual of the homology class represented by the
map $s_{\alpha _{i}}$, $1\leq i\leq n$, then one has that

\bigskip

\noindent \textbf{Lemma 2.4. }\textsl{The set }$\left\{ \omega _{1},\cdots
,\omega _{n}\right\} $\textsl{\ is a basis of the cohomology group }$%
H^{2}(G/T)$.$\square $

\bigskip

On the other hand let $\Theta =\{\theta _{1},\cdots ,\theta _{n}\}$ be an
ordered basis of the unit lattice $\Lambda _{e}$. It defines $n$ oriented
circle subgroups on the maximal torus

\begin{enumerate}
\item[(2.2)] $\widetilde{\theta }_{i}:S^{1}=\mathbb{R}/\mathbb{Z}\rightarrow
T$ by $\widetilde{\theta }_{i}(t):=\exp (t\theta _{i})$, $1\leq i\leq n$,
\end{enumerate}

\noindent that represent also an ordered basis of the first homology $%
H_{1}(T)$. As result if we let $t_{i}\in H^{1}(T)$ be the class Kronnecker
dual to the map $\widetilde{\theta }_{i}$, then

\begin{enumerate}
\item[(2.3)] $H^{\ast }(T)=\Lambda (t_{1},\cdots ,t_{n})$ (i.e. the exterior
ring generated by $t_{1},\cdots ,t_{n}$).
\end{enumerate}

In view of the relation $\Lambda _{r}\subseteq \Lambda _{e}$ by Lemma 2.3
there exists an unique integer matrix $C(\Theta )=\left( c_{i,j}\right)
_{n\times n}$ expressing the ordered basis $\Delta $ of $\Lambda _{r}$ by
the ordered basis $\Theta $ of $\Lambda _{e}$. That is, the relation

\begin{quote}
$\qquad \left( 
\begin{tabular}{l}
$\alpha _{1}$ \\ 
$\vdots $ \\ 
$\alpha _{n}$%
\end{tabular}%
\right) =C(\Theta )\left( 
\begin{tabular}{l}
$\theta _{1}$ \\ 
$\vdots $ \\ 
$\theta _{n}$%
\end{tabular}%
\right) $
\end{quote}

\noindent holds on $L(T)$. Our main result is

\bigskip

\noindent \textbf{Theorem 2.5.}\textsl{\ With respect to the basis\ }$%
\left\{ t_{1},\cdots ,t_{n}\right\} $ \textsl{on} $H^{1}(T)$\textsl{\ and
and the basis }$\left\{ \omega _{1},\cdots ,\omega _{n}\right\} $\textsl{\
on }$H^{2}(G/T)$\textsl{, the transgression }$\tau $ \textsl{is given by the
formula}

\begin{enumerate}
\item[(2.4)] $\left( 
\begin{tabular}{l}
$\tau (t_{1})$ \\ 
$\vdots $ \\ 
$\tau (t_{n})$%
\end{tabular}%
\right) =C(\Theta )^{\tau }\left( 
\begin{tabular}{l}
$\omega _{1}$ \\ 
$\vdots $ \\ 
$\omega _{n}$%
\end{tabular}%
\right) $\textsl{,}
\end{enumerate}

\noindent \textsl{where }$C(\Theta )^{\tau }$\textsl{\ is the transpose of
the matrix }$C(\Theta )$\textsl{.}

\bigskip

\noindent \textbf{Proof. }We begin with the simple case where the group $G$
is $1$--connected. Then a basis $\Theta $ of the unit lattice $\Lambda
_{e}=\Lambda _{r}$ can be taken to be $\Delta =\{\alpha _{1},\cdots ,\alpha
_{n}\}$ by i) of Lemma 2.3. Since $C(\Theta )$ is now the identity matrix we
are bound to show that

\begin{quote}
$\tau (t_{i})=\omega _{i}$, $1\leq i\leq n$.
\end{quote}

For a simple root $\alpha _{i}\in $ $\Delta $ the inclusion $(K(\alpha
_{i}),S^{1})\subset (G,T)$ of subgroups induces the following bundle map
over $s_{\alpha _{i}}$:

\begin{center}
\begin{tabular}{lll}
$\quad S^{1}$ & $\overset{\widetilde{\alpha }_{i}}{\rightarrow }$ & $T$ \\ 
$\quad \downarrow $ &  & $\downarrow $ \\ 
$K(\alpha _{i})=S^{3}$ & $\rightarrow $ & $G$ \\ 
$\pi _{i}\downarrow $ &  & $\downarrow \pi $ \\ 
$K(\alpha _{i})/S^{1}=S^{2}$ & $\overset{s_{\alpha _{i}}}{\rightarrow }$ & $%
G/T$%
\end{tabular}%
,
\end{center}

\noindent where, since the group $G$ is $1$--connected, the subgroup $%
K(\alpha _{i})$ is isomorphic to the $3$--sphere $S^{3}$ (i.e. the group of
unit quaternions), while the map $\pi _{i}$ is the Hopf fibration over $%
S^{2} $. This implies that, in the homotopy exact sequence of the firation $%
\pi $

\begin{quote}
$\cdots \rightarrow \pi _{2}(G)\overset{\pi _{\ast }}{\rightarrow }\pi
_{2}(G/T)\overset{\partial }{\rightarrow }\pi _{1}(T)\rightarrow \cdots $,
\end{quote}

\noindent the connecting homomorphism $\partial $ satisfies the relation

\begin{enumerate}
\item[(2.5)] $\partial \left[ s_{\alpha _{i}}\right] =\left[ \widetilde{%
\alpha }_{i}\right] $, $1\leq i\leq n$.
\end{enumerate}

\noindent Since both of the Hurewicz homomorphisms

\begin{quote}
$\pi _{1}(T)\rightarrow H_{1}(T)$ and $\pi _{2}(G/T)\rightarrow H_{2}(G/T)$
\end{quote}

\noindent are isomorphisms, and since the transgression $\tau $ is Kronecker
dual to $\partial $ in the sense that

\begin{quote}
$\tau =Hom(\partial ,1):Hom(\pi _{1}(T),\mathbb{Z})\rightarrow Hom(\pi
_{2}(G/T),\mathbb{Z})$
\end{quote}

\noindent one obtains from (2.5) that $\tau (t_{i})=$ $\omega _{i}$, $1\leq
i\leq n$.

Turning to a general situation assume that the group $G$ is semisimple. Let $%
d:(G_{0},T_{0})\rightarrow (G,T)$ be the universal cover of $G$ with $T_{0}$
the maximal torus on $G_{0}$ that corresponds to $T$ under $d$. Then, with
respect to the canonical identifications (induced by the tangent map of $d$
at the group unit)

\begin{quote}
$L(G_{0})=L(G)$ and $L(T_{0})=L(T)$,
\end{quote}

\noindent the exponential map of $G$ admits the decomposition

\begin{quote}
$\exp =d\circ \exp _{0}:(L(G_{0}),L(T_{0}))\rightarrow
(G_{0},T_{0})\rightarrow (G,T)$,
\end{quote}

\noindent where $\exp _{0}$ is the exponential map of $G_{0}$. It follows
that, if we let

\begin{quote}
$p(\Lambda _{r},\Lambda _{e}):$ $T_{0}=L(T_{0})/\Lambda _{r}$ $\rightarrow
T=L(T_{0})/\Lambda _{e}$
\end{quote}

\noindent be the covering map induced by the inclusion $\Lambda
_{r}\subseteq \Lambda _{e}$ of the lattices, then

\begin{enumerate}
\item[(2.6)] $d\mid T_{0}=p(\Lambda _{r},\Lambda _{e}):T_{0}\rightarrow T$.
\end{enumerate}

\noindent Note that the induced map $p(\Lambda _{r},\Lambda _{e})_{\ast }$
on $\pi _{1}(T_{0})$ is determined by the matrix $C(\Theta )=\left(
c_{ij}\right) _{n\times n}$ as

\begin{enumerate}
\item[(2.7)] $p(\Lambda _{r},\Lambda _{e})_{\ast }[\widetilde{\alpha }%
_{i}]=c_{i,1}[\widetilde{\theta }_{1}]+\cdots +c_{i,n}[\widetilde{\theta }%
_{n}]$.
\end{enumerate}

On the other hand, by the naturality of homotopy exact sequence of
firations, the restriction $d\mid T_{0}$ fits into the commutative diagram

\begin{enumerate}
\item[(2.8)] $%
\begin{array}{ccc}
\pi _{2}(G_{0}/T_{0}) & \underset{\cong }{\overset{\partial _{0}}{%
\rightarrow }} & \pi _{1}(T_{0}) \\ 
\parallel &  & \qquad \downarrow (d\mid T_{0})_{\ast } \\ 
\pi _{2}(G/T) & \overset{\partial }{\rightarrow } & \pi _{1}(T)%
\end{array}%
$,
\end{enumerate}

\noindent where the vertical identification on the left comes from the fact
that the covering $d:(G_{0},T_{0})\rightarrow (G,T)$ induces a
diffeomorphism $G_{0}/T_{0}\cong G/T$, and where $\partial _{0}$, $\partial $
are the connecting homomorphisms in the homotopy exact sequences of the
bundles $G_{0}\rightarrow G_{0}/T_{0}$, $G\rightarrow G/T$, respectively. It
follows that, for a simple root $\alpha _{i}\in \Delta $, one has

\begin{quote}
$\partial \left[ s_{\alpha _{i}}\right] =(d\mid T_{0})_{\ast }\circ \partial
_{0}\left[ s_{\alpha _{i}}\right] $ (by the diagram (2.8))

$\qquad =(d\mid T_{0})_{\ast }\left[ \widetilde{\alpha }_{i}\right] $ (by
the proof of the previous case)

$\qquad =p(\Lambda _{r},\Lambda _{e})_{\ast }(\left[ \widetilde{\alpha }_{i}%
\right] )$ (by (2.6)).
\end{quote}

\noindent The proof is now completed by (2.7), together with the fact
(again) that the map $\tau $ is Kronecker dual to $\partial $.$\square $

\section{Applications}

In a concrete situation formula (2.4) is ready to apply to evaluate the
transgression $\tau $ (henceforth, the differential $d_{2}$ on $E_{2}^{\ast
,\ast }(G)$) associated to the fibration $G\rightarrow G/T$. We present
below three examples.

Assume firstly that the group $G$ is $1$--connected. By i) of Lemma 2.3 one
can take the set $\Delta =\{\alpha _{1},\cdots ,\alpha _{n}\}$ of simple
roots as a preferable basis of the unit lattice $\Lambda _{e}$. The
transition matrix $C(\Theta )$ is then the identity matrix. Theorem 2.5
implies that

\bigskip

\noindent \textbf{Corollary 3.1. }\textsl{If the group }$G$\textsl{\ is }$1$%
\textsl{--connected, there exists a basis }$\left\{ t_{1},\cdots
,t_{n}\right\} $\textsl{\ on }$H^{1}(T)$\textsl{\ so that }$\tau (t_{i})=$%
\textsl{\ }$\omega _{i}$\textsl{, }$1\leq i\leq n$\textsl{.}$\square $

\bigskip

Suppose next that the group $G$ is adjoint. According to ii) of Lemma 2.3
the set $\Omega =\{\phi _{1},\cdots ,\phi _{n}\}$ of fundamental dominant
weights is a basis of $\Lambda _{e}$, and the transition matrix $C(\Theta )$
from $\Lambda _{e}$ to $\Lambda _{r}$ is\textsl{\ }the Cartan matrix $A$ of $%
G$ by iii) of Lemma 2.3. We get from Theorem 2.5 that

$\bigskip $

\noindent \textbf{Corollary 3.2. }\textsl{If the group }$G$\textsl{\ is
adjoint, there exists a basis }$\left\{ t_{1},\cdots ,t_{n}\right\} $\textsl{%
\ on }$H^{1}(T)$\textsl{\ so that}

\begin{quote}
$\left( 
\begin{tabular}{l}
$\tau (t_{1})$ \\ 
$\vdots $ \\ 
$\tau (t_{n})$%
\end{tabular}%
\right) =A^{\tau }\left( 
\begin{tabular}{l}
$\omega _{1}$ \\ 
$\vdots $ \\ 
$\omega _{n}$%
\end{tabular}%
\right) $\textsl{,}
\end{quote}

\noindent \textsl{where }$A$\textsl{\ is the Cartan matrix of the group }$G$%
\textsl{.}$\square $

\bigskip

Since both of the groups $H^{1}(T)$ and $H^{2}(G/T)$ are torsion free,
formula (2.4) is also applicable to evaluate the transgression

\begin{quote}
$\tau :H^{1}(T;\mathbb{Z}_{p})\rightarrow H^{2}(G/T;\mathbb{Z}_{p})$
\end{quote}

\noindent for the cohomologies with coefficients in the cyclic ring $\mathbb{%
Z}_{p}$ of order $p>1$. Given a Lie group $G$ denote by $PG:=G/\mathcal{Z}%
(G) $ the associated Lie group of the adjoint type, where $\mathcal{Z}(G)$
is the center of $G$. In what follows we shall assume that $G$ is one of the
the $1$--connected Lie groups $Sp(n)$, $E_{6}$ or $E_{7}$, and shall adopt
the following conventions:

\begin{quote}
i) a set $\Delta =\{\alpha _{1},\cdots ,\alpha _{n}\}$ of simple roots of $G$
is given and ordered as the vertex of the Dykin diagram of $G$ pictured on 
\cite[p.58]{H};

ii) the preferable basis of the unit lattice $\Lambda _{e}$ of the group $PG$
is taken to be $\Omega $.
\end{quote}

\noindent Granted with the Cartan matrix of $G$ presented on \cite[p.59]{H}
Corollary 3.2 implies that

\bigskip

\noindent \textbf{Corollary 3.3. }\textsl{For each of the pairs }$(G,p)=$%
\textsl{\ }$(Sp(n),2)$\textsl{, }$(E_{6},3)$\textsl{\ and }$(E_{7},2)$%
\textsl{, the kernel and cokernel of the transgression}

\begin{quote}
$\tau :H^{1}(T;\mathbb{Z}_{p})\rightarrow H^{2}(PG/T;\mathbb{Z}_{p})$,
\end{quote}

\noindent \textsl{are both isomorphic to }$\mathbb{Z}_{p}$\textsl{, whose
generators are specified in the following table}

\begin{quote}
\begin{tabular}{l||l|l}
\hline
$(G,p)$ & Generator of $\ker \tau =\mathbb{Z}_{p}$ & Generator of co$\ker
\tau =\mathbb{Z}_{p}$ \\ \hline\hline
$(Sp(n),2)$ & $t_{n}$ & $\omega _{1}$ \\ \hline
$(E_{6},3)$ & $t_{1}-t_{3}+t_{5}-t_{6}$ & $\omega _{1}$ \\ \hline
$(E_{7},2)$ & $t_{2}+t_{5}+t_{7}$ & $\omega _{2}$ \\ \hline
\end{tabular}%
\textsl{.}$\square $
\end{quote}

\bigskip

\noindent \textbf{Remark 3.4. }Let\textbf{\ }$G$ be a compact semisimple Lie
group and let $p$ be prime. In \cite[formula (4)]{K} Ka\v{c} stated a
formula for the differential $d_{2}$ on $E_{2}^{\ast ,\ast }(G;\mathbb{Z}%
_{p})$ which implies that the transgression $\tau $ in the characteristic $p$
is always an isomorphism. This contradicts to Corollary 3.3.

On the other hand, in the context of Schubert calculus a unified
presentation of the cohomology ring of the quotient manifolds $G/T$ has been
obtained in \cite[Theorem 1.2]{DZ1}. Combining this result with Theorem 2.5
of the present paper, method to construct the integral cohomology rings of
compact connected Lie groups using the spectral sequence $\left\{
E_{r}^{\ast ,\ast }(G),d_{r}\right\} $ has been developed in the works \cite%
{D,DZ2}.$\square $

\bigskip

\textbf{Acknowledgement.} The author is grateful to the referee for valuable
suggestions.

\end{document}